\documentclass[12pt]{amsart}

\usepackage{amssymb,euscript,amsmath, mathrsfs,amscd,mathtools}
\usepackage[dvips]{graphicx}
\usepackage{epsfig,pstricks}
\usepackage{fullpage}

\newcounter{ENUM}
\newcommand{\be}{\begin{enumerate}}
\newcommand{\ee}{\end{enumerate}}
\newcommand{\beas}{\begin{eqnarray*}}
\newcommand{\eeas}{\end{eqnarray*}}
\newcommand{\bea}{\begin{eqnarray}}
\newcommand{\eea}{\end{eqnarray}}
\newcommand{\beq}{\begin{equation}}
\newcommand{\eeq}{\end{equation}}

\newcommand{\st}{\,:\,}

\newcommand{\sn}{\mathfrak{S}_n}
\newcommand{\fs}{\mathfrak{S}}

\newcommand{\caln}{\mathcal{N}}

\newtheorem{thm}{Theorem}[section]

\newtheorem{conj}[thm]{Conjecture}

\theoremstyle{definition}

\theoremstyle{remark}

\numberwithin{equation}{section}

\def\zz{\mathbb{Z}}
\def\cc{\mathbb{C}}

\newcommand{\qq}{\mathbb{Q}}


\begin{document}
\title{Some Schubert Shenanigans}

\date{\today}

\author{Richard P. Stanley}
\email{rstan@math.mit.edu}
\address{Department of Mathematics, University of Miami, Coral Gables,
FL 33124}

%\thanks{Partially supported by NSF grant DMS-1068625.}

\begin{abstract}
We give a conjectured evaluation of the determinant of a certain
matrix $\tilde{D}(n,k)$. The entries of $\tilde{D}(n,k)$ are
either 0 or specializations $\fs_w(1,\dots,1)$ of Schubert
polynomials. The conjecture implies that the weak order of the
symmetric group $S_n$ has the strong Sperner property. A number of
peripheral results and problems are also discussed.
\end{abstract}

\maketitle

\section{Introduction} \label{intro}
The primary purpose of this paper is to present a conjectured
evaluation of the determinant of a certain matrix $\tilde{D}(n,k)$,
where $0\leq k< \frac 12\binom n2$, related to the weak order $W_n$ of
the symmetric group $S_n$. The entries of $\tilde{D}(n,k)$ are either
0 or $\fs_w(1,\dots,1)$, the sum of the coefficients of the Schubert
polynomial $\fs_w$, for certain permutations $w$.  If
$\det\tilde{D}(n,k)\neq 0$ for all $0\leq k<\frac 12\binom n2$, then
$W_n$ has the strong Sperner property. It is currently unknown whether
$W_n$ has just the ordinary Sperner property. We have an explicit
conjecture for the value $\det\tilde{D}(n,k)$ which shows in particular
that it is nonzero, but we have no idea how to prove this conjecture.

The above problem suggests some other questions related to the numbers
  \beq \nu_w\coloneqq \fs_w(1,1,\dots,1). \label{eq:nu} \eeq
In particular, it is well-known that $\nu_w=1$ if and only if
$w$ is 132-avoiding (or \emph{dominant}). We have a conjectured
condition for when $\nu_w=2$ and can prove the sufficiency of this
condition. One can go on to consider when $\nu_w=3$, etc., though a
``nice'' characterization for every $m$ of when $\nu_w=m$ does not
seem very likely.

The consideration above of small values of $\nu_w$ suggests looking at
how large $\nu_w$ can be for $w\in S_n$. Write $u(n)$ for this maximum
value. We have some rather crude bounds
which show that as $n\to\infty$,
 $$ \frac 14\leq\liminf \frac{\log_2 u(n)}{n^2}\leq 
   \limsup \frac{\log_2 u(n)}{n^2}\leq \frac 12. $$
In addition to improving these bounds, one can also ask if there is
some kind of ``limiting shape'' for the permutation(s) $w\in S_n$
satisfying $\nu_w=u(n)$.

\section{A weak order determinant}
For enumeration and poset terminology see \cite{ec1}. For any set $X$,
we denote by $\qq X$ the $\qq$-vector space with basis $X$. Let
$P=P_0\cup P_1\cup\cdots \cup P_m$ be a finite
graded rank-symmetric poset of rank $m$. Thus if $p_k=\#P_k$ then
$p_k=p_{m-k}$. A linear transformation $U\colon \qq P\to \qq P$ is
\emph{order-raising} if for every $t\in P$ we have $U(t)\in \qq
C^+(t)$, where $C^+(t)$ denotes the set of elements of $P$ that cover
$t$. Thus in this case
$U$ maps $\qq P_k$ into $\qq P_{k+1}$. The poset $P$ is
\emph{rank-unimodal} if
  $$ p_0\leq p_1\leq\cdots\leq p_j\geq p_{j+1}\geq\cdots\geq p_m $$
for some $j$. Note that if also $P$ is rank-symmetric, then we can
take $j=\lfloor m/2\rfloor$. We say that $P$ is \emph{strongly
  Sperner} if for all $r\geq 1$, the largest subset $S$ of $P$
containing no $(r+1)$-element chain has the same number of elements as
the largest union $T$ of $r$ levels of $P$ (in which case we can take
$S=T$). The following result is given in \cite[Lemma~1.1]{rs:hl}.

\begin{thm} \label{thm:sperner}
Suppose there exists an order-raising operator
$U\colon \qq P\to \qq P$ such that if $0\leq k<\frac m2$
then the linear transformation $U^{m-2k}\colon \qq P_k\to \qq P_{m-k}$
is a bijection. Then $P$ is strongly Sperner.
\end{thm}

We would like to apply Theorem~\ref{thm:sperner} to the weak order
$W_n$ of the symmetric group $S_n$. A permutation $w$ has rank $k$ in
$W_n$, denoted $w\in (W_n)_k$, if $w$ has $k$ inversions. The
rank-generating function of $W_n$ is given by
 $$ F(W_n,q)\coloneqq \sum_{k=0}^{\binom n2}\#(W_n)_k\cdot q^k
      = (1+q)(1+q+q^2)\cdots (1+q+\cdots+q^{n-1}). $$
Let $\ell(w)$ denote the number of inversions (or \emph{length}) of
$w$. Write $s_i$ for the adjacent
transposition $(i,i+1)$, so $v$ covers $u$ in $W_n$ if and only if
$v=us_i$ for some $1\leq i\leq n-1$ for which $\ell(v)= 1+\ell(u)$. We
now define the order-raising operator $U$ by
 $$ U(u) =\sum_{i\st \ell(us_i)=1+\ell(u)} i\cdot us_i,\ \ u\in
      W_n. $$ 
It follows from \cite[(6.11)]{macd} or \cite[Lemma~2.3]{f-s} that
for $u\in W_n$ and $j\geq 0$ we have
  \beq U^j(u) = \left(\binom n2-2k\right)!
   \sum_v \nu_{vu^{-1}}v, \label{eq:dnkuv} \eeq
where $v$ ranges over all elements in $W_n$ satisying $\ell(v) =
\ell(u) +j$ and $v>u$ (in weak order), and $\nu_{vu^{-1}}$ is defined
by equation~\eqref{eq:nu}.

Let $D(n,k)$ denote the matrix of the linear transformation
$U^{\binom n2-2k}\colon \qq (W_n)_k\to \qq (W_n)_{\binom n2-k}$ with respect
to the bases $(W_n)_k$ and $(W_n)_{\binom n2-k}$ (in some order). The
idea of looking at this matrix is implicit in \cite{rs:mo}. By
equation~\eqref{eq:dnkuv} the $(u,v)$-entry of this matrix (for $u\in
(W_n)_k$ and $v\in (W_n)_{\binom n2-k}$) is given by
   $$ D(n,k)_{uv}=\left\{ \begin{array}{rl}
     \left(\binom n2-2k\right)!\nu_{vu^{-1}}, & \mathrm{if}\ u\leq
     v\ \mathrm{in}\ W_n\\[.5em] 
     0, & \mathrm{otherwise}. \end{array} \right. $$

We can divide each entry of $D(n,k)$ by $\left(\binom n2-2k\right)!$,
obtaining a matrix $\tilde{D}(n,k)$ satisfying
  $$ \det D(n,k) = \left(\binom
       n2-2k\right)!^{^{^{\scriptsize\#(W_n)_k}}} \det 
    \tilde{D}(n,k). $$
The following conjecture has been verified for all pairs $(n,k)$
satisfying both $n\leq 12$ and $k\leq 5$, as well as a few other
cases.

\begin{conj} \label{conj:det}
  We have
  \beq \det \tilde{D}(n,k) = \pm\prod_{i=0}^{k-1}\left( \frac{\binom
    n2-(k+i)} {k-i}\right)^{\#(W_n)_i}. \label{eq:conj} \eeq
\end{conj}

We can prove Conjecture~\ref{conj:det} when $k=1$. This amounts to
showing that
 \beq \det\left[ \begin{array}{ccccccccc}
    1 & 1 & 1 & \cdots & 1 & 1 & 2 & 0\\
    1 & 1 & 1 & \cdots & 1 & 2 & 0 & 1\\
    1 & 1 & 1 & \cdots & 2 & 0 & 1 & 1\\
      &   &   & \vdots\\
    1 & 1 & 2 & \cdots & 1 & 1 & 1 & 1\\
    1 & 2 & 0 & \cdots & 1 & 1 & 1 & 1\\
    2 & 0 & 1 & \cdots & 1 & 1 & 1 & 1\\
    0 & 1 & 1 & \cdots & 1 & 1 & 1 & 1
  \end{array} \right] = (-1)^{\lfloor (n-2)/2\rfloor}
     \left(\binom n2-1\right), \label{eq:k1det} \eeq
where the matrix is $(n-1)\times (n-1)$.

Conjecture~\ref{conj:det} can be formulated in terms of the
\emph{nilCoxeter algebra} $\caln_n$ of the symmetric group $S_n$. By
definition, $\caln_n$ is a $\qq$-algebra with basis $\{\sigma_w\st
w\in S_n\}$ and relations 
$$ \sigma_u\sigma_v = \left\{ \begin{array}{rl}
  \sigma_{uv}, & \mathrm{if}\ \ell(uv)=\ell(u)+\ell(v)\\[.5em]
   0, & \mathrm{otherwise}. \end{array} \right. $$
The algebra $\caln_n$ is graded via
 $$ \caln_n = (\caln_n)_0\oplus (\caln_n)_1\oplus\cdots\oplus
   (\caln_n)_{\binom n2}, $$
where $(\caln_n)_i$ is spanned by elements $\sigma_w$ for which
$\ell(w)=i$. Thus
 $$ (\caln_n)_i(\caln_n)_j\subseteq (\caln_n)_{i+j}. $$
Define
  $$ \theta=\sum_{i=1}^{n-1}i\sigma_{s_i}\in (\caln_n)_1, $$
where $s_i$ denotes the adjacent transposition $(i,i+1)$.
Then for $0\leq k<\frac 12\binom n2$, the matrix of the map
$\qq(\caln_n)_k\to \qq(\caln_n)_{\binom n2-k}$, with respect to the
bases $(\caln_n)_k$ and $(\caln_n)_{\binom n2-k}$, given by
multiplication by $\theta^{\binom n2-2k}$ coincides with
$D(n,k)$. Thus Conjecture~\ref{conj:det} implies that $\theta$ is a
kind of ``hard Lefschetz element'' for $\caln_n$. Note, however, that
the commutativity properties of $\caln_n$ imply that it is not the
cohomology ring of a topological space.

The right-hand side of equation~\eqref{eq:conj} is somewhat
reminiscent of formulas arising in the theory of differential posets
\cite[Sect.~4]{rs:dp}
\cite[pp.~154--157]{rs:vdp}. This leads us to ask
whether there might be some ``down operator'' $D\colon \qq W_n
\to \qq W_n$ taking $\qq (W_n)_{k+1}$ to $\qq (W_n)_k$ with ``nice''
commutativity properties with respect to 
$U$.  Ideally the linear transformation $DU\colon \qq W_n\to
\qq W_n$ should have integer eigenvalues. We have been unable to find
a suitable candidate for $D$. For the weak order of \emph{affine}
Coxeter groups  there is a nice theory due to Lam and Shimozono
\cite{l-s}. 

\section{Generalizations and refinements}
We mentioned that the entries of the matrix $\tilde{D}(n,k)$ are
either 0 or specializations $\nu_w$ of Schubert polynomials
$\fs_w$. It is natural to ask what happens when we replace $\nu_w$
with $\fs_w$ itself. However, we have been unable to discover anything
interesting in this regard.

We can also consider a $q$-analogue of
Conjecture~\ref{conj:det}. Theorem~2.4 of \cite{f-s} (originally
conjectured by Macdonald in \cite[(6.11$_q$?)]{macd}) suggests that
the correct $q$-analogue of $\tilde{D}(n,k)$ is obtained by replacing
the entry $\fs_w(1,1,\dots,1)$ of $\tilde{D}(n,k)$ with
$\fs_w(1,q,q^2,\dots,q^{m-1})$ (where $\fs_w$ is a polynomial in
$x_1,\dots, x_m$). However, the data is not too  encouraging. For instance,
when $(n,k)=(5,2)$ the determinant is $q^{36}$ times an irreducible
polynomial (over $\qq$) of degree 20. For $(n,k)=(5,3)$ the
determinant has the form
  $$ \pm q^{26}(q^2+q+1)^3(q^4+q^3+q^2+q+1)^3 P(q), $$
where $P(q)$ is an irreducible polynomial of degree 28. Perhaps
$\fs_w(1,q,q^2,\dots,q^{m-1})$ needs to be multiplied by some power of $q$.

It is natural to ask what happens when we replace the weak order $W_n$
with the (strong) Bruhat order on $S_n$, which we continue to denote
as $S_n$. For the Bruhat order the strong Sperner property is quite easy to
see \cite[p.~182]{rs:hl}. However, it is still interesting to ask for
an analogue $E(n,k)$ of the matrix $D(n,k)$. There are reasons to
define the operator $V\colon \qq S_n\to \qq S_n$ as follows. Write
$t_{ij}$ for the permutation in $S_n$ that transposes $i$ and
$j$. Then for $u\in S_n$ define
  $$ V(u) = \sum_{\substack{1\leq i<j\leq n\\ \ell(ut_{ij})=1+\ell(u)}}
      (j-i)\cdot ut_{ij}. $$
One reason for this definition is the fact, due essentially to
Chevalley and further investigated by Stembridge \cite{stemb} and
Postnikov and Stanley \cite{p-s}, that
   $$ V^{\binom n2}(\mathrm{id}) = \binom n2!\,w_0, $$
where id denotes the identity permutation and $w_0=n,n-1,\dots,1$.
Moreover, for $u\in (S_n)_j$ and $v\in (S_n)_{j+1}$ we have
 \beq [v]V(u) = [\fs_v]\fs_u(\fs_{s_1}+\cdots+\fs_{s_{n-1}}),
 \label{eq:vvu} \eeq
where $[v]V(u)$ denotes the coefficient of $v$ in $V(u)$ and
$[\fs_v]F$ denotes the coefficient of $\fs_v$ when the polynomial $F$
is written as a linear combination of Schubert polynomials. As above
$s_i$ denotes the adjacent transposition $(i,i+1)$, so
$\fs_{s_i}=x_1+x_2+\cdots +x_i$.
   
Let $E(n,k)$ denote the matrix of the linear transformation
$V^{n-2k}\colon \qq (S_n)_k\to \qq (S_n)_{\binom n2-k}$ with respect
to the bases $(S_n)_k$ and $(S_n)_{\binom n2-k}$ (in some order). Let
$e(n,k) =\det E(n,k)$. For $u\in (S_n)_k$ and $v\in (S_n)_{\binom
  n2-2k}$, it follows from equation~\eqref{eq:vvu} that the
$(u,v)$-entry of $E(n,k)$ is given by 
  $$ E(n,k)_{uv} = [\fs_v]
     \fs_u(\fs_{s_1}+\cdots+\fs_{s_{n-1}})^{\binom n2-2k}. $$
We have computed that
  \beas e(4,1) & = & \pm 2^7\cdot 3\cdot 5^2\cdot 19\\
        e(4,2) & = & \pm 2^6\cdot 3\cdot 29\\
        e(5,1) & = & \pm 2^{22}\cdot 3^6\cdot 5^5 \cdot 7^4
                     \cdot 59\cdot 89. \eeas
Obviously some further computation, or even better a theorem, is needed.

Unlike the situation for the weak order, the Bruhat analogue
$\mathcal{B}_n$ of $\mathcal{N}_n$ is the cohomology ring of a
(smooth, projective) variety, viz., the complete flag variety
$\mathcal{F}_n = \mathrm{GL}(n,\cc)/B$, where $B$ is a Borel
subgroup. Thus the hard Lefschetz theorem for $\mathcal{F}_n$
immediately implies the strong Sperner property for $S_n$
\cite{rs:hl}, but as stated above this is easy to prove directly. A
variety related to $\mathcal{F}_n$ is  the Grassmann variety Gr$(n,k)$
of $k$-dimensional subspaces of $\cc^n$. Here the determinants
of the matrices analogous to $D(n,k)$ and $E(n,k)$ have been computed
by Proctor \cite{proctor}. He also computed determinants corresponding
to products of complex projective spaces, where the corresponding
poset is a product of chains.

A canonical question suggested by Conjecture~\ref{conj:det} is its
possible extension to other finite Coxeter groups. We have not made
any computations in this direction.

Another possible extension of Conjecture~\ref{conj:det} is to refine
the determinant of $\tilde{D}(n,k)$. The two most natural refinements
are the characteristic polynomial (or equivalently the eigenvalues) of
$\tilde{D}(n,k)$ and the Smith normal form (SNF) over $\zz$ of
$\tilde{D}(n,k)$. The characteristic polynomial depends on the order
of the rows and columns. We have done some computation using every
possible order and have found nothing interesting. For instance, the
eigenvalues are in general not integers. On the other hand, the SNF
does not depend on the order of the rows and columns, and the data
looks intriguing. For instance, the ratio between two consecutive
diagonal terms is always quite small. Write $f(n,k)$ for the vector of
diagonal entries of the SNF of $\tilde{D}(n,k)$. Here is some data
(which we have not attempted to analyze, except for $k=1$ discussed
below). We use exponents to denote repetition, e.g., $(1^5,3^2)$
stands for $(1,1,1,1,1,3,3)$.

   \beas f(4,1) & = & (1^2,5)\\
   f(4,2) & = & (1^2,3^2,6)\\
   f(5,1) & = & (1^3,9)\\
   f(5,2) & = & (1^5,7^3,28)\\
   f(5,3) & = & (1^6,5^6,15^2,105)\\
   f(5,4) & = & (1^8,3^4,6^4,30^4)\\
   f(6,1) & = & (1^4,14)\\
   f(6,2) & = & (1^9,6,12^3,156)\\
   f(6,3) & = & (1^{15},5^4,10^5,110^4,220)\\
   f(6,4) & = & (1^{20},2,4^9,8^5,24^5,72^4,360^4,3960)\\
   f(6,5) & = & (1^{31},3^6,6^5,42^{19},84,168^4,504^5)\\
   f(6,6) & = & (1^{28},2^{18},10^2,20^{28},140^{10},280^3,840)\\
   f(6,7) & = & (1^{52},2^{18},6^{10},12^6,60^{11},420^3,840)\\
   f(7,1) & = & (1^5,20)\\
   f(7,2) & = & (1^{14},9,18^4,342)\\
   f(7,3) & = & (1^{29},8^5,16^9,272^5,816)\\
   f(7,4) & = & (1^{49},7^{14},14^{15},70^6,210^8,420,1680^4,28560)\\
   f(7,5) & = &
    (1^{76},2^9,6^{20},12^{15},156^{29},1092^{15},5460^4,21840)\\
   f(8,1) & = & (1^6,27)\\
   f(8,2) & = & (1^{20},25^6,325)\\
   f(8,3) & = & (1^{49},23^{20},92,276^5,6900)\\
   f(8,4) & = & (1^{98},7^6,21^{43},231^{20},5313^6,10626)\\
   f(9,1) & = & (1^7,35)\\
   f(9,2) & = & (1^{27},33^7,561)\\
   f(9,3) & = & (1^{76},31^{27},496^7,5456)\\
   f(10,1) & = & (1^8,44)\\
   f(10,2) & = & (1^{35},21,42^7,1806)
   \eeas
The data suggests the obvious conjecture $f(n,1)=(1^{n-2},\binom
n2-1)$. We can prove this as follows. We can take $\tilde{D}(n,1)$ to
be the matrix of equation~\eqref{eq:k1det}. Since $\det
\tilde{D}(n,1)=\binom n2-1$ (see equation~\eqref{eq:k1det}), by basic
properties of SNF (e.g., \cite[Thm.~2.4]{rs:snf}) it suffices to show
that if $\tilde{D}'(n,1)$ denotes $\tilde{D}(n,1)$ with its first row
and column deleted, then $\det \tilde{D}'(n,1)=1$. The proof is by
induction on $n$, the base case $n=3$ being trivial to check. Subtract
the next-to-last row of $\tilde{D}'(n,1)$ from the last row and expand
by the last row. We obtain $\det \tilde{D}'(n,1)=\det
\tilde{D}'(n-1,1)$, so the proof follows.

   \section{Two-term Schubert polynomials}
The appearance of the numbers $\nu_w=\fs_w(1,1,\dots,1)$ as entries of
the matrix $\tilde{D}(n,k)$ suggest the question of how ``complicated'' these
numbers can be. The simplest situation is when $\fs_w$ is a single
monomial, so $\nu_w=1$. It is well-known \cite{macd} that $\nu_w=1$ if
and only if  $w$ is 132-avoiding (or \emph{dominant}).

\begin{conj} \label{conj:two}
We have $\nu_w=2$ if and only if $w=a_1\cdots a_n$ has exactly one
subsequence $a_i,a_j,a_k$ in the pattern 132, i.e., $a_i<a_k<a_j$. 
\end{conj}

If $f(n)$ denotes the number of permutations with exactly one
subsequence with pattern 132, then it is known that
$f(n)=\binom{2n-3}{n-3}$, $n\geq 3$ \cite[eqn.~(2)]{bona}.

We can prove the ``if'' direction of
Conjecture~\ref{conj:two}. Namely, if $w=a_1\cdots a_n$ has exactly
one subsequence with pattern 132, then it is easily seen that the
unique such subsequence has the form $a_i,a_{i+1},a_j$. If we
transpose $a_i$ and 
$a_{i+1}$, then we get a permutation $w'$ that is 132-avoiding. For
any permutation $v=v_1\cdots v_n$ let
  $$ \lambda_i(v) = \#\{j>i\st v_j<v_i\}. $$
By the theory of Schubert
polynomials \cite{macd}, we have (since $w'$ is 132-avoiding)
  \beq \fs_{w'} = x_1^{\lambda_1(w')}x_2^{\lambda_2(w')}\cdots
  x_{n-1}^{\lambda_{n-1}(w')} \label{eq:dom} \eeq
and (since $w$ covers $w'$ in $W_n$)
  \beq \fs_w = \frac{\fs_{w'}(x_1,\dots,x_i,x_{i+1},\dots,x_{n-1})-
    \fs_{w'}(x_1,\dots,x_{i+1},x_i,\dots,x_{n-1})}{x_i-x_{i+1}}.
     \label{eq:fsw} \eeq
Since $a_i,a_{i+1},a_j$ is the unique subsequence of $w$ with the
pattern 132, we have $\lambda_i(w')=\lambda_{i+1}(w')+2$. It follows
from equations~\eqref{eq:dom} and \eqref{eq:fsw} that $\fs_w$ has two
terms. 
    
We suspect that the converse, i.e., if $\nu_w=2$ then $w$ has exactly
one subsequence with pattern 132, will not be difficult to prove. 

One can go on to consider the number $f_j(n)$ of permutations
$w\in\sn$ for which $\nu_w=j\geq 3$, and the characterization of
such permutations. It seems unlikely that there will be a nice result
in general.  

\section{The maximum value of $\nu_w$ for $w\in S_n$}
The previous section dealt with small values of $\nu_w$, so we can
also ask for large values. As in Section~\ref{intro}, define
  $$ u(n) =\max \{ \nu_w\st w\in S_n\}. $$
Setting each $x_i=y_i=1$ in the ``Cauchy identity''
 $$ \prod_{i+j\leq n} (x_i+y_j) = \sum_{\substack{w=v^{-1}u\\
    \ell(w)=\ell(u)+\ell(v)}} \fs_u(x)\fs_v(y) $$
immediately yields the estimates
   $$ \frac 14\leq\liminf \frac{\log_2 u(n)}{n^2}\leq 
   \limsup \frac{\log_2 u(n)}{n^2}\leq \frac 12. $$
Presumably the limit $L=\lim_{n\to\infty}n^{-2}\log_2 u(n)$ exists (in
which case $\frac 14\leq L\leq \frac 12$), and it would be interesting to
compute this limit. It would also be interesting to determine whether
the permutations $w\in S_n$ achieving $\nu_e=u(n)$ have some
``limiting'' description. Below is the value of $u(n)$ for $3\leq
n\leq 9$, together with the set of all permutations $w\in S_n$
achieving $\nu_w=u(n)$.

\begin{center}
 \begin{tabular}{crc}
   $n$ & $u(n)$ & $w$\\ \hline
   3 & 2 & 132\\
   4 & 5 & 1432\\
   5 & 14 & 12543,\ 15432,\ 21543\\
   6 & 84 & 126543,\ 216543\\
   7 & 660 & 1327654\\
   8 & 9438 & 13287654\\
   9 & 163592 & 132987654\\
   10 & 4424420 & 1,4,3,2,10,9,8,7,6,5
 \end{tabular}
\end{center}

\end{document}